\documentclass[12pt]{article}

\catcode`\@=11

\thicklines
\newskip\Einheit \Einheit=.6cm
\newcount\xcoord \newcount\ycoord
\newdimen\xdim \newdimen\ydim \newdimen\PfadD@cke \newdimen\Pfadd@cke
\def\PfadDicke#1{\PfadD@cke#1 \divide\PfadD@cke by2 
\Pfadd@cke\PfadD@cke \multiply\PfadD@cke by2}
\long\def\LOOP#1\REPEAT{\def\BODY{#1}\ITERATE}
\def\ITERATE{\BODY \let\next\ITERATE \else\let\next\relax\fi \next}
\let\REPEAT=\fi
\def\Punkt{\hbox{\raise-2pt\hbox to0pt{\hss\scriptsize$\bullet$\hss}}}

\def\DuennPunkt(#1,#2){\unskip
  \raise#2 \Einheit\hbox to0pt{\hskip#1 \Einheit
          \raise-1.5pt\hbox to0pt{\hss\tiny$\bullet$\hss}\hss}}
		  
\def\NormalPunkt(#1,#2){\unskip
  \raise#2 \Einheit\hbox to0pt{\hskip#1 \Einheit
          \raise-3pt\hbox to0pt{\hss\large$\bullet$\hss}\hss}}
\def\DickPunkt(#1,#2){\unskip
  \raise#2 \Einheit\hbox to0pt{\hskip#1 \Einheit
          \raise-4pt\hbox to0pt{\hss\Large$\bullet$\hss}\hss}}
\def\Kreis(#1,#2){\unskip
  \raise#2 \Einheit\hbox to0pt{\hskip#1 \Einheit
          \raise-4pt\hbox to0pt{\hss\Large$\circ$\hss}\hss}}
\def\Diagonale(#1,#2)#3{\unskip\leavevmode
  \xcoord#1\relax \ycoord#2\relax
      \raise\ycoord \Einheit\hbox to0pt{\hskip\xcoord \Einheit
         \unitlength\Einheit
         \line(1,1){#3}\hss}}
\def\AntiDiagonale(#1,#2)#3{\unskip\leavevmode
  \xcoord#1\relax \ycoord#2\relax \advance\xcoord by -0.05\relax
      \raise\ycoord \Einheit\hbox to0pt{\hskip\xcoord \Einheit
         \unitlength\Einheit
         \line(1,-1){#3}\hss}}
\def\Pfad(#1,#2),#3\endPfad{\unskip\leavevmode
  \xcoord#1 \ycoord#2 \thicklines\ZeichnePfad#3\endPfad\thinlines}
\def\ZeichnePfad#1{\ifx#1\endPfad\let\next\relax
  \else\let\next\ZeichnePfad
    \ifnum#1=1
      \raise\ycoord \Einheit\hbox to0pt{\hskip\xcoord \Einheit
         \vrule height\Pfadd@cke width1 \Einheit depth\Pfadd@cke\hss}%
      \advance\xcoord by 1
    \else\ifnum#1=2
      \raise\ycoord \Einheit\hbox to0pt{\hskip\xcoord \Einheit
        \hbox{\hskip-\PfadD@cke\vrule height1 \Einheit 
width\PfadD@cke depth0pt}\hss}%
      \advance\ycoord by 1
    \else\ifnum#1=3
      \raise\ycoord \Einheit\hbox to0pt{\hskip\xcoord \Einheit
         \unitlength\Einheit
         \line(1,1){1}\hss}
      \advance\xcoord by 1
      \advance\ycoord by 1
    \else\ifnum#1=4
      \raise\ycoord \Einheit\hbox to0pt{\hskip\xcoord \Einheit
         \unitlength\Einheit
         \line(1,-1){1}\hss}
      \advance\xcoord by 1
      \advance\ycoord by -1
    \fi\fi\fi\fi
  \fi\next}
\def\hSSchritt{\leavevmode\raise-.4pt\hbox 
to0pt{\hss.\hss}\hskip.2\Einheit
  \raise-.4pt\hbox to0pt{\hss.\hss}\hskip.2\Einheit
  \raise-.4pt\hbox to0pt{\hss.\hss}\hskip.2\Einheit
  \raise-.4pt\hbox to0pt{\hss.\hss}\hskip.2\Einheit
  \raise-.4pt\hbox to0pt{\hss.\hss}\hskip.2\Einheit}
\def\vSSchritt{\vbox{\baselineskip.2\Einheit\lineskiplimit0pt
\hbox{.}\hbox{.}\hbox{.}\hbox{.}\hbox{.}}}
\def\DSSchritt{\leavevmode\raise-.4pt\hbox to0pt{%
  \hbox to0pt{\hss.\hss}\hskip.2\Einheit
  \raise.2\Einheit\hbox to0pt{\hss.\hss}\hskip.2\Einheit
  \raise.4\Einheit\hbox to0pt{\hss.\hss}\hskip.2\Einheit
  \raise.6\Einheit\hbox to0pt{\hss.\hss}\hskip.2\Einheit
  \raise.8\Einheit\hbox to0pt{\hss.\hss}\hss}}
\def\dSSchritt{\leavevmode\raise-.4pt\hbox to0pt{%
  \hbox to0pt{\hss.\hss}\hskip.2\Einheit
  \raise-.2\Einheit\hbox to0pt{\hss.\hss}\hskip.2\Einheit
  \raise-.4\Einheit\hbox to0pt{\hss.\hss}\hskip.2\Einheit
  \raise-.6\Einheit\hbox to0pt{\hss.\hss}\hskip.2\Einheit
  \raise-.8\Einheit\hbox to0pt{\hss.\hss}\hss}}
\def\SPfad(#1,#2),#3\endSPfad{\unskip\leavevmode
  \xcoord#1 \ycoord#2 \ZeichneSPfad#3\endSPfad}
\def\ZeichneSPfad#1{\ifx#1\endSPfad\let\next\relax
  \else\let\next\ZeichneSPfad
    \ifnum#1=1
      \raise\ycoord \Einheit\hbox to0pt{\hskip\xcoord \Einheit
         \hSSchritt\hss}%
      \advance\xcoord by 1
    \else\ifnum#1=2
      \raise\ycoord \Einheit\hbox to0pt{\hskip\xcoord \Einheit
        \hbox{\hskip-2pt \vSSchritt}\hss}%
      \advance\ycoord by 1
    \else\ifnum#1=3
      \raise\ycoord \Einheit\hbox to0pt{\hskip\xcoord \Einheit
         \DSSchritt\hss}
      \advance\xcoord by 1
      \advance\ycoord by 1
    \else\ifnum#1=4
      \raise\ycoord \Einheit\hbox to0pt{\hskip\xcoord \Einheit
         \dSSchritt\hss}
      \advance\xcoord by 1
      \advance\ycoord by -1
    \fi\fi\fi\fi
  \fi\next}
\def\Koordinatenachsen(#1,#2){\unskip
 \hbox to0pt{\hskip-.5pt\vrule height#2 \Einheit width.5pt depth1 
\Einheit}%
 \hbox to0pt{\hskip-1 \Einheit \xcoord#1 \advance\xcoord by1
    \vrule height0.25pt width\xcoord \Einheit depth0.25pt\hss}}
\def\Koordinatenachsen(#1,#2)(#3,#4){\unskip
 \hbox to0pt{\hskip-.5pt \ycoord-#4 \advance\ycoord by1
    \vrule height#2 \Einheit width.5pt depth\ycoord \Einheit}%
 \hbox to0pt{\hskip-1 \Einheit \hskip#3\Einheit 
    \xcoord#1 \advance\xcoord by1 \advance\xcoord by-#3 
    \vrule height0.25pt width\xcoord \Einheit depth0.25pt\hss}}
\def\Gitter(#1,#2){\unskip \xcoord0 \ycoord0 \leavevmode
  \LOOP\ifnum\ycoord<#2
    \loop\ifnum\xcoord<#1
      \raise\ycoord \Einheit\hbox to0pt{\hskip\xcoord 
\Einheit\Punkt\hss}%
      \advance\xcoord by1
    \repeat
    \xcoord0
    \advance\ycoord by1
  \REPEAT}
\def\Gitter(#1,#2)(#3,#4){\unskip \xcoord#3 \ycoord#4 \leavevmode
  \LOOP\ifnum\ycoord<#2
    \loop\ifnum\xcoord<#1
      \raise\ycoord \Einheit\hbox to0pt{\hskip\xcoord 
\Einheit\Punkt\hss}%
      \advance\xcoord by1
    \repeat
    \xcoord#3
    \advance\ycoord by1
  \REPEAT}
\def\Label#1#2(#3,#4){\unskip \xdim#3 \Einheit \ydim#4 \Einheit
  \def\lo{\advance\xdim by-.5 \Einheit \advance\ydim by.5 \Einheit}%
  \def\llo{\advance\xdim by-.25cm \advance\ydim by.5 \Einheit}%
  \def\loo{\advance\xdim by-.5 \Einheit \advance\ydim by.25cm}%
  \def\o{\advance\ydim by.25cm}%
  \def\ro{\advance\xdim by.5 \Einheit \advance\ydim by.5 \Einheit}%
  \def\rro{\advance\xdim by.25cm \advance\ydim by.5 \Einheit}%
  \def\roo{\advance\xdim by.5 \Einheit \advance\ydim by.25cm}%
  \def\l{\advance\xdim by-.30cm}%
  \def\r{\advance\xdim by.30cm}%
  \def\lu{\advance\xdim by-.5 \Einheit \advance\ydim by-.6 \Einheit}%
  \def\llu{\advance\xdim by-.25cm \advance\ydim by-.6 \Einheit}%
  \def\luu{\advance\xdim by-.5 \Einheit \advance\ydim by-.30cm}%
  \def\u{\advance\ydim by-.30cm}%
  \def\ru{\advance\xdim by.5 \Einheit \advance\ydim by-.6 \Einheit}%
  \def\rru{\advance\xdim by.25cm \advance\ydim by-.6 \Einheit}%
  \def\ruu{\advance\xdim by.5 \Einheit \advance\ydim by-.30cm}%
  #1\raise\ydim\hbox to0pt{\hskip\xdim
     \vbox to0pt{\vss\hbox to0pt{\hss$#2$\hss}\vss}\hss}%
}
\catcode`\@=12

\usepackage{amsmath,amsthm}
\usepackage{color}
\usepackage{xspace}
\usepackage[colorlinks=true,
linkcolor=green,
filecolor=brown,
citecolor=green]{hyperref}
\def\red{\textcolor{red} }
\def\green{\textcolor{green} }

\begin{document}
\newtheorem{lemma}{Lemma}
\newtheorem{theorem}{Theorem}
\newtheorem{prop}{Proposition}
\newtheorem{cor}{Corollary}
\begin{center}
{\Large
Two Bijections for Dyck Path Parameters \\ 
}
\vspace{10mm}
DAVID CALLAN  \\
Department of Statistics  \\
University of Wisconsin-Madison  \\
1210 W. Dayton St   \\
Madison, WI \ 53706-1693  \\
{\bf callan@stat.wisc.edu}  \\
\vspace{5mm}
July 19 2004
\end{center}

\vspace{5mm}

The Motzkin number $M_{n}$
(\htmladdnormallink{A001006}{http://www.research.att.com:80/cgi-bin/access.cgi/as/njas/sequences/eisA.cgi?Anum=A001006})
counts Motzkin $n$-paths: lattice paths of upsteps $U=(1,1)$, 
flatsteps $F=(1,0)$ and downsteps $D=(1,-1)$ such that (i) the path 
contains $n$ steps, (ii) the number of $U$s =  number of $D$s, and 
(iii) the path never dips below the 
horizontal line joining its initial and terminal points (ground level) \cite{ec2}. 
A Dyck $n$-path is a Motzkin $(2n)$-path with no flatsteps, counted 
by the Catalan number $C_{n}$. A 
$UUU$-free  Dyck path is one that contains no consecutive $UUU$ and 
so on. A descent in a path is a maximal sequence of 
contiguous downsteps. A short descent is one consisting of a 
single downstep. A terminal descent is one that ends the path. Thus a Dyck path is $UDU$-free iff it contains no 
short nonterminal descents. Each upstep in a Dyck (or Motzkin) path has 
a matching downstep: the first one encountered directly East of the 
upstep.

The first result is easy, and also appears in \cite{ElizaldeMansour03}.

\begin{theorem}
    The number of $UUU$-free  Dyck $n$-paths is $M_{n}$.
\end{theorem}

\begin{proof}
    Given a $UUU$-free  Dyck $n$-path, change each 
$UUD$ to $U$, then change each remaining $UD$ to $F$. This is a bijection to 
Motzkin $n$-paths. For example with $n=5: U\,U\,D\,U\,D\,U\,U\,D\,D\,D -> 
U\,F\,U\,D\,D$. The inverse is obvious.
\end{proof}

The next bijection gives the distributions both of the parameter ``number 
of $UDU$s'' 
(\htmladdnormallink{A091869}{http://www.research.att.com:80/cgi-bin/access.cgi/as/njas/sequences/eisA.cgi?Anum=A091869})
\cite{deutschMotzkin} and the parameter ``number 
of $DDU$s'' 
(\htmladdnormallink{A091894}{http://www.research.att.com:80/cgi-bin/access.cgi/as/njas/sequences/eisA.cgi?Anum=A091894}) 
on Dyck paths.
\newpage

\begin{theorem} \ \\ 
\vspace*{-5mm}

     $($i$\,)$   The number of Dyck $n$-paths containing exactly $k$ $UDU$s is 
    $\binom{n-1}{k}M_{n-1-k}$\\ $($Donaghey distribution$)$.
    
 $($ii$)$   The number of Dyck $n$-paths containing exactly $k$ $DDU$s is 
    $\binom{n-1}{2k}2^{n-1-2k}C_{k}$\\ $($Touchard distribution$)$.
\end{theorem}

\begin{proof}
A bicolored Motzkin $n$-path is a Motzkin path of length $n$ in which 
each flatstep is given one of two colors, say green and black, denoted 
by wavy and flat overlines respectively. There is an obvious bijection 
to Dyck $(n+1)$-paths: replace each $U$ by $UU$, $D$ by $DD$, 
$\overline{F}$ by $UD$, $\widetilde{F}$ by $UD$ and then prepend a $U$ 
and append a $D$. 

Here is a less obvious bijection that takes $\#\,\widetilde{F}$s to 
$\#\,UDU$s and $\#\,D$s to $\#\,DDU$s. Since it is easy to count 
bicolored Motzkin paths either by $\#\,\widetilde{F}$s or by $\#\,D$s, we 
get the stated distributions.

Given a bicolored Motzkin $n$-path, first append a downstep as a convenience. 
Now every flatstep $F$ has an 
associated downstep: the first $D$ whose initial point is on the ray 
obtained by extending $F$ eastward. 
Then 
\vspace*{-2mm}
\begin{enumerate}
    \item  Leave $U$s intact

    \item  Replace $D$ by $UDD$

    \item  Replace $\widetilde{F}$ by $UD$

    \item  Replace $\overline{F}$ by $U$ and insert a $D$ immediately 
    before its associated downstep (thus the new $U$ and inserted $D$ 
    will be a matching $UD$ pair and the result is the same no matter 
    in what order 
    the $\overline{F}$s are processed). 
\end{enumerate} \vspace*{-2mm}
Finally, delete the appended $D$ (it remains undisturbed). An example with $n=14$ is illustrated.
\vspace*{2mm}
\Einheit=0.5cm
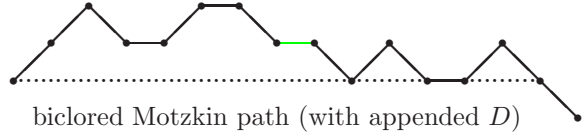
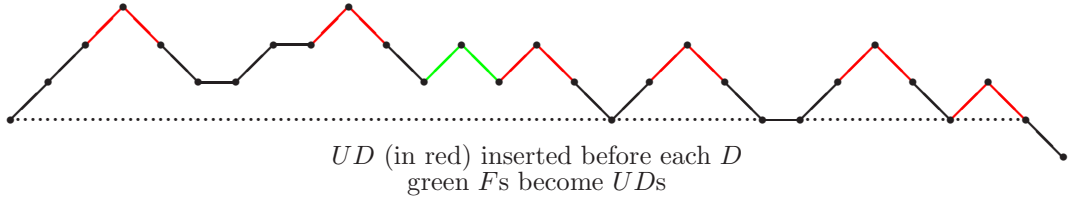
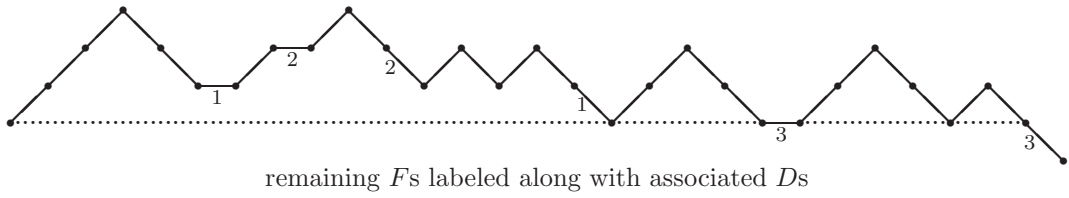
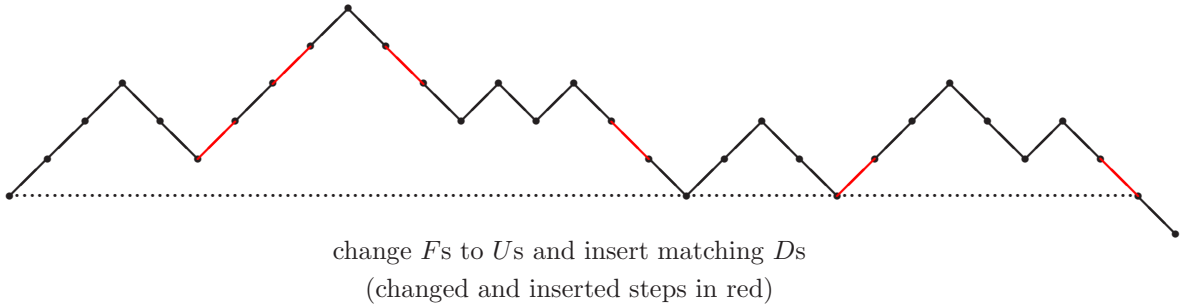
\[
\Pfad(-7,7),3341314\endPfad
\green{\Pfad(0,8),1\endPfad}
\Pfad(1,8),4341344\endPfad
\SPfad(-7,7),11111111111\endSPfad
\SPfad(5,7),11\endSPfad
\DuennPunkt(-7,7)
\DuennPunkt(-6,8)
\DuennPunkt(-5,9)
\DuennPunkt(-4,8)
\DuennPunkt(-3,8)
\DuennPunkt(-2,9)
\DuennPunkt(-1,9)
\DuennPunkt(0,8)
\DuennPunkt(1,8)
\DuennPunkt(2,7)
\DuennPunkt(3,8)
\DuennPunkt(4,7)
\DuennPunkt(5,7)
\DuennPunkt(6,8)
\DuennPunkt(7,7)
\DuennPunkt(8,6)
\Label\o{\textrm{\footnotesize{biclored Motzkin path (with appended $D$)}}}(0,5.5)
\Pfad(-14,0),33\endPfad
\Pfad(-10,2),4131\endPfad
\Pfad(-4,2),4\endPfad
\Pfad(1,1),43\endPfad
\Pfad(5,1),413\endPfad
\Pfad(10,1),4\endPfad
\Pfad(13,0),4\endPfad
\green{\Pfad(-3,1),34\endPfad}
\red{
\Pfad(-12,2),34\endPfad
\Pfad(-6,2),34\endPfad
\Pfad(-1,1),34\endPfad
\Pfad(3,1),34\endPfad
\Pfad(8,1),34\endPfad
\Pfad(11,0),34\endPfad
}
\SPfad(-14,0),11111111111111111111\endSPfad
\SPfad(7,0),111111\endSPfad
\DuennPunkt(-14,0)
\DuennPunkt(-13,1)
\DuennPunkt(-12,2)
\DuennPunkt(-11,3)
\DuennPunkt(-10,2)
\DuennPunkt(-9,1)
\DuennPunkt(-8,1)
\DuennPunkt(-7,2)
\DuennPunkt(-6,2)
\DuennPunkt(-5,3)
\DuennPunkt(-4,2)
\DuennPunkt(-3,1)
\DuennPunkt(-2,2)
\DuennPunkt(-1,1)
\DuennPunkt(0,2)
\DuennPunkt(1,1)
\DuennPunkt(2,0)
\DuennPunkt(3,1)
\DuennPunkt(4,2)
\DuennPunkt(5,1)
\DuennPunkt(6,0)
\DuennPunkt(7,0)
\DuennPunkt(8,1)
\DuennPunkt(9,2)
\DuennPunkt(10,1)
\DuennPunkt(11,0)
\DuennPunkt(12,1)
\DuennPunkt(13,0)
\DuennPunkt(14,-1)
\Label\o{\textrm{\footnotesize{$UD$ (in red) inserted before each $D$}}}(0,-1.5)
\Label\o{\textrm{\footnotesize{green $F$s become $UD$s}}}(0,-2.2)
\] 

\vspace*{10mm}

\Einheit=0.5cm
\[
\Pfad(-14,0),3334413134434\endPfad    
\Pfad(-1,1),344334413344344\endPfad
\SPfad(-14,0),11111111111111111111\endSPfad
\SPfad(7,0),111111\endSPfad
\DuennPunkt(-14,0)
\DuennPunkt(-13,1)
\DuennPunkt(-12,2)
\DuennPunkt(-11,3)
\DuennPunkt(-10,2)
\DuennPunkt(-9,1)
\DuennPunkt(-8,1)
\DuennPunkt(-7,2)
\DuennPunkt(-6,2)
\DuennPunkt(-5,3)
\DuennPunkt(-4,2)
\DuennPunkt(-3,1)
\DuennPunkt(-2,2)
\DuennPunkt(-1,1)
\DuennPunkt(0,2)
\DuennPunkt(1,1)
\DuennPunkt(2,0)
\DuennPunkt(3,1)
\DuennPunkt(4,2)
\DuennPunkt(5,1)
\DuennPunkt(6,0)
\DuennPunkt(7,0)
\DuennPunkt(8,1)
\DuennPunkt(9,2)
\DuennPunkt(10,1)
\DuennPunkt(11,0)
\DuennPunkt(12,1)
\DuennPunkt(13,0)
\DuennPunkt(14,-1)
\Label\o{\textrm{\footnotesize{remaining $F$s labeled along with associated 
$D$s}}}(0,-2)
\Label\o{\textrm{\scriptsize{1}}}(-8.5,0.2)
\Label\o{\textrm{\scriptsize{2}}}(-6.5,1.2)
\Label\o{\textrm{\scriptsize{3}}}(6.5,-.8)
\Label\o{\textrm{\scriptsize{1}}}(1.2,0)
\Label\o{\textrm{\scriptsize{2}}}(-3.9,1)
\Label\o{\textrm{\scriptsize{3}}}(13.1,-1)
\] 

\vspace*{10mm}
\Einheit=0.5cm
\[
\Pfad(-15,0),33344\endPfad
\Pfad(-9,2),3\endPfad
\Pfad(-7,4),34\endPfad
\Pfad(-4,3),434\endPfad
\Pfad(-1,2),34\endPfad
\Pfad(2,1),43344\endPfad
\Pfad(8,1),334434\endPfad
\Pfad(15,0),4\endPfad
\SPfad(-15,0),111111111111111111111111111111\endSPfad
\DuennPunkt(-15,0)
\DuennPunkt(-14,1)
\DuennPunkt(-13,2)
\DuennPunkt(-12,3)
\DuennPunkt(-11,2)
\DuennPunkt(-10,1)
\DuennPunkt(-9,2)
\DuennPunkt(-8,3)
\DuennPunkt(-7,4)
\DuennPunkt(-6,5)
\DuennPunkt(-5,4)
\DuennPunkt(-4,3)
\DuennPunkt(-3,2)
\DuennPunkt(-2,3)
\DuennPunkt(-1,2)
\DuennPunkt(0,3)
\DuennPunkt(1,2)
\DuennPunkt(2,1)
\DuennPunkt(3,0)
\DuennPunkt(4,1)
\DuennPunkt(5,2)
\DuennPunkt(6,1)
\DuennPunkt(7,0)
\DuennPunkt(8,1)
\DuennPunkt(9,2)
\DuennPunkt(10,3)
\DuennPunkt(11,2)
\DuennPunkt(12,1)
\DuennPunkt(13,2)
\DuennPunkt(14,1)
\DuennPunkt(15,0)
\DuennPunkt(16,-1)
\red{
\Pfad(-10,1),3\endPfad
\Pfad(-8,3),3\endPfad
\Pfad(7,0),3\endPfad
\Pfad(-5,4),4\endPfad
\Pfad(1,2),4\endPfad
\Pfad(14,1),4\endPfad
}
\Label\o{\textrm{\footnotesize{change $F$s to $U$s and insert matching $D$s }}}(0,-2)
\Label\o{\textrm{\footnotesize{(changed and inserted steps in red)  }}}(0,-3)
\] 

The map is invertible: $\widetilde{F}$s are recaptured as $UD$s 
followed by a $U$, $D$s are recaptured as $UDD$s, $\overline{F}$s are 
recaptured as $U$s whose matching $D$ is in the \emph{interior} of a 
descent of length $\ge 3$ and the original $U$s are recaptured as $U$s 
such that both the $U$ itself and its matching $D$ are followed by a $U$.
\end{proof}

This map restricted to Motzkin paths (no green $F$s) is a bijection 
from Motzkin $n$-paths to $UDU$-free Dyck $(n+1)$-paths.
Further, note that a Motzkin path contains no flatsteps at ground level iff 
the corresponding $UDU$-free Dyck path ends with $UD$. So, by deleting this $UD$, Motzkin 
$n$-paths with no flatsteps at ground level
(\htmladdnormallink{A005043}{http://www.research.att.com:80/cgi-bin/access.cgi/as/njas/sequences/eisA.cgi?Anum= A005043})
correspond to Dyck $n$-paths with no short descents.

\end{document}